\documentclass[12pt]{amsart}
\usepackage{amscd,amssymb}
\newtheorem{theorem}{Theorem}
\newtheorem{lemma}[theorem]{Lemma}
\newtheorem{proposition}[theorem]{Proposition}
\newtheorem{corollary}[theorem]{Corollary}
\theoremstyle{definition}

\newtheorem{problem}[theorem]{Problem}
\newtheorem{example}[theorem]{Example}

\theoremstyle{remark}
\newtheorem{remark}[theorem]{Remark}

\def\varph{{\varphi}}
\def\ra{{\rightarrow}}
\def\lra{{\longrightarrow}}
\def\om{{\omega}}

\def\la{{\lambda}}
\def\al{{\alpha}}

\def\ga{{\gamma}}

\def\vGa{{\varGamma}}

\def\bZ{{\mathbb Z}}
\def\bQ{{\mathbb Q}}
\def\bR{{\mathbb R}}
\def\inv{{^{-1}}}
\def\Sg{{\Sigma_g}}
\def\Mg{{\mathcal M}_g}
\def\M{{\mathcal M}}
\def\Msg{{\mathcal M}_{g,*}}
\newcommand\Hom{\operatorname{Hom}}

\newcommand\Flux{\operatorname{Flux}}
\newcommand\Symp{\operatorname{Symp}}
\newcommand\Ham{\operatorname{Ham}}

\newcommand\eFlux{\operatorname{\widetilde{Flux}}}
\newcommand\CAL{\operatorname{Cal}}
\newcommand\Ker{\operatorname{Ker}}

\newcommand\Int{\operatorname{Int}}
\newcommand\sign{\operatorname{sign}}
\newcommand\BS{\operatorname{BSymp}}
\newcommand\ES{\operatorname{ESymp}}
\newcommand\BD{\operatorname{BDiff}}
\newcommand\ED{\operatorname{EDiff}}
\newcommand\BED{\operatorname{BEDiff}}
\newcommand\BG{\operatorname{B\bar\Gamma}}
\newcommand\Diff{\operatorname{Diff}}

\newcommand\id{\operatorname{Id}}
\catcode`\@=\active 

\begin{document}

\title[Signatures of foliated surface bundles]
{Signatures of foliated surface bundles and
the symplectomorphism groups of surfaces}
\author{D.~Kotschick}
\address{Mathematisches Institut, Ludwig-Maximilians-Universit\"at M\"unchen,
Theresienstr.~39, 80333 M\"unchen, Germany}
\email{dieter{\char'100}member.ams.org}
\author{S.~Morita}
\address{Department of Mathematical Sciences\\
University of Tokyo \\Komaba, Tokyo 153-8914\\
Japan}
\email{morita{\char'100}ms.u-tokyo.ac.jp}


\thanks{Support from the {\sl Deutscher Akademischer Austauschdienst}
and the {\sl Deutsche Forschungsgemeinschaft} is gratefully acknowledged.
The first named author is a member of the {\sl European Differential
Geometry Endeavour} (EDGE), Research Training Network HPRN-CT-2000-00101,
supported by The European Human Potential Programme. The second named author is
partially supported by JSPS Grant 13440017}
\date{May 8, 2003; MSC 2000 classification: primary 57R17, 57R30, 57R50;
secondary 37E30, 57M99, 58H10}

\begin{abstract}
For any closed oriented surface $\Sigma_g$ of genus $g\geq 3$, we
prove the existence of {\it foliated} $\Sigma_g$-bundles over surfaces
such that the signatures of the total spaces are non-zero. We can
arrange that the total holonomy of the horizontal foliations
preserve a prescribed symplectic form
$\om$ on the fiber. We relate the cohomology class represented by the
transverse symplectic form to a {\it crossed} homomorphism
$\eFlux:\Symp \Sigma_g\ra H^1(\Sigma_g;\bR)$ which is an
extension of the flux homomorphism
$\Flux:\Symp_0 \Sigma_g\ra H^1(\Sigma_g;\bR)$
from the identity component $\Symp_0\Sg$ to the whole group
$\Symp \Sigma_g$ of symplectomorphisms of $\Sigma_g$ with respect to
the symplectic form $\om$.
\end{abstract}

\maketitle

\section{Statement of the main results}\label{s:intro}

Let $\Sg$ be a closed oriented surface of genus $g$.
It is a classical result that, for any $g\geq 3$, there exist
oriented $\Sg$-bundles over closed oriented surfaces such that
the signatures of the total spaces are non-zero, see
Kodaira~\cite{Kodaira} and Atiyah~\cite{Atiyah}.
In this paper, we prove the existence of such bundles which,
in addition to having non-zero signature, are flat, or foliated.
This means that there exist codimension two foliations
complementary to the fibers, which is equivalent to the existence
of lifts of the holonomy homomorphisms from the mapping class
group to the diffeomorphism group of the fiber.
We will further show that such lifts
can be chosen to preserve a prescribed area form,
or equivalently a symplectic form $\om$, on the fiber.
More precisely, we prove the following result.
\begin{theorem}\label{th:main}
For any $g\geq 3$, there exist foliated oriented $\Sg$-bundles
$\pi : E \ra B$ over closed oriented surfaces $B$
such that the total holonomy group is contained in
the symplectomorphism group $\Symp \Sg$ 
with respect to a prescribed symplectic
form $\om$ on $\Sg$, and $\sign E \not= 0$.
\end{theorem}
Our proof in Section~\ref{s:proof1} below is not constructive at the
final stage.
In particular, we do not have any explicit example of a foliated
$\Sg$-bundle with non-zero signature.
Also, we know no examples of surface bundles over surfaces that can 
be shown not to admit any foliated structure. Thurston's
result that Haefliger's classifying space $\BG_2$ is $3$-connected,
see~\cite{Thurston2,Thurston}, implies that the tangent bundle along the
fibers of any surface bundle over a surface is homotopic to the
normal bundle of a codimension two foliation on the total space.
However, it is unclear whether we can arrange this
foliation to be transverse to the fibers everywhere.

Let $\pi:E\ra B$ be a foliated oriented $\Sg$-bundle
as in Theorem~\ref{th:main}, so that
$\sign E\not= 0$ and the image of the total holonomy homomorphism
$$
\pi_1 B\lra \Diff_{+}\Sg
$$
is contained in the symplectomorphism subgroup $\Symp\Sg\subset
\Diff_{+}\Sg$ with respect to $\om$. Since the total holonomy
preserves the symplectic form $\om$ on $\Sg$, the pullback of this
form to the product $\Sg\times\tilde B$ descends to
$E=(\Sg\times\tilde B)/\pi_{1}B$ as a globally defined closed $2$-form
$\tilde{\om}$ of rank $2$ which restricts to $\om$ on the fiber.
Hence we have the corresponding cohomology class
$$
v=[\tilde{\om}]\in H^{2}(E;\bR) \ ,
$$
which we call the {\it transverse symplectic class}.
At the universal space level, this cohomology class $v$
can be considered as an element of
$
H^2(\ES^{\delta}\Sg;\bR)$,
where the discrete group $\ES^{\delta}\Sg$ is defined
as follows.
Let $\Mg$ and $\Msg$ denote the mapping class group of $\Sg$, 
respectively the mapping class group relative to a base point. 
Then we have the universal extension $\pi_1\Sg\ra\Msg\ra\Mg$. 
If we pull back this extension by the natural projection
$\Symp^{\delta}\Sg\ra\Mg$, where the symbol $\delta$ indicates
the {\it discrete} topology, we obtain an extension
$$
1\lra\pi_1\Sg\lra\ES^{\delta}\Sg\lra\Symp^{\delta}\Sg\lra 1 \ .
$$
Thus, $\ES^{\delta}\Sg$ is the universal model for the
fundamental groups of the total spaces of
foliated $\Sg$-bundles with area-preserving holonomy.

On the other hand, we have the Euler class $e\in H^{2}(E;\bZ)$
of the tangent bundle along the fibers of $\pi$. This bundle is
the normal bundle of the horizontal foliation on $E$. The two
cohomology classes $v$ and $e$ are proportional on the fiber,
and if we normalize $\om$ so that
$$
\int_{\Sg}\om=2g-2 \ ,
$$
then $e+v$ restricts to $0$ on the fiber. However,
we can never have the equality $v=-e$ for
of the following reason. Clearly, we have
$v^2=0$ (since $\tilde{\om}^2$ vanishes identically),
while $e^2\not= 0$ since its fiber integral is nothing but
the first Mumford--Morita--Miller class
$e_1\in H^2(X;\bZ)$ which represents the signature of
the total space (see~\cite{Atiyah,Meyer,Morita2}). Thus the
question of identifying the difference of $v$ and $-e$ arises.
We shall answer this by making use of certain basic facts
in symplectic topology.

Let $\Symp_0\Sg$ denote the identity component of $\Symp\Sg$.
Then there is a well-defined surjective homomorphism
$$
\Flux: \Symp_0\Sg\lra H^1(\Sg;\bR) \ ,
$$
called the flux homomorphism.
We refer to the book~\cite{MS} by McDuff and Salamon for generalities of
symplectic topology, including the flux homomorphism
as well as the Calabi homomorphism used in the proof of
Theorem~\ref{th:main3} below.

\begin{theorem}\label{th:main2}
For all $g\geq 2$, the flux homomorphism
can be extended to a crossed homomorphism
$$
\widetilde{\Flux}: \Symp\Sg\lra H^1(\Sg;\bR) \ .
$$
This extension is unique in the sense that its cohomology class in
$
H^1(\Symp^{\delta}\Sg;H^1(\Sg;\bR))
$
is unique.
Furthermore, in the cohomology spectral sequence of the
extension
$$
1\lra\pi_1\Sg\lra\ES^{\delta}\Sg\lra\Symp^{\delta}\Sg\lra 1
$$
the class
$$
e+v \in \Ker\big(H^2(\ES^{\delta}\Sg;\bR)\lra H^2(\Sg;\bR)\big)
$$
projects to the above cohomology class
$$
[\widetilde{\Flux}]\in H^1(\Symp^{\delta}\Sg;H^1(\Sg;\bR)) \ .
$$
\end{theorem}

We will actually determine the group $H^1(\Symp^{\delta}\Sg;H^1(\Sg;\bR))$
completely, see Proposition~\ref{prop:h1}. Furthermore, in Section~\ref{s:final}
we generalize Theorem~\ref{th:main2} to a certain class of closed 
symplectic manifolds of higher dimensions, see Proposition~\ref{prop:fluxhi}.

Theorem~\ref{th:main} can be reformulated and extended
in the context of the Mumford--Morita--Miller classes
(see~\cite{Mumford,Morita1,Miller}).
Let $\Mg$ be the mapping class group of $\Sg$ as before and
let $e_i=\pi_{*}e^{i+1}\in H^{2i}(\Mg;\bQ)$ be the $i^{\text{th}}$
Mumford--Morita--Miller class with {\it rational}
coefficients defined by integration over the fiber in the universal
$\Sg$-bundle.

Let $\Diff_{+}\Sg$ be the group of orientation-preserving diffeomorphisms
of $\Sg$. Then $\Mg$ can be considered as the group of path components of
$\Diff_{+}\Sg$ and we have an extension
$$
1\lra\Diff_0\Sg\lra\Diff_{+}\Sg\overset{p}{\lra}\Mg\lra 1 \ ,
$$
where $\Diff_0\Sg$ is the identity component of
$\Diff_{+}\Sg$ and $p$ is the natural projection.
It follows from the Bott vanishing theorem for the characteristic
classes of the normal bundles of foliations that
$$
p^*(e_i)=0\in H^{2i}(\BD^{\delta}_{+}\Sg;\bQ)
$$
for all $i\geq 3$, see~\cite{Morita2} and also~\cite{Morita4}.
(Here, as before, $\delta$ indicates the {\it discrete}
topology, so that the space $\BD^{\delta}_{+}\Sg$ is the
classifying space of {\it foliated} oriented $\Sg$-bundles.) More
precisely, the Bott vanishing theorem applied to the horizontal
foliation shows that $e^{i}=0\in
H^{2i}(\BED^{\delta}_{+}\Sg;\bQ)$ for all $i\geq 4$, where
$\ED^{\delta}_{+}\Sg$ is defined as in the case of symplectomorphism
groups considered above. On the other hand,
Theorem~\ref{th:main} shows that $e^{2}\in
H^4(\BED^{\delta}_{+}\Sg;\bQ)$ and its fiber integral
$e_{1}\in H^2(\BD^{\delta}_{+}\Sg;\bQ)$ are non-zero.

It remains to determine whether other polynomials in
$e_1$ and $e_2$ are trivial in $H^*(\BD^{\delta}_{+}\Sg;\bQ)$,
or not. By extending Theorem~\ref{th:main}, we can give a
partial answer to this question. Namely, we show the
non-triviality of any power
$e_1^k\in H^{2k}(\BD^{\delta}_{+}\Sg;\bQ)$ of the
first characteristic class $e_1$.
In fact, we can prove the following stronger
non-vanishing result for the subgroup
$\Symp\Sg\subset\Diff_{+}\Sg$.

\begin{theorem}\label{th:main3}
Let $\Symp^{\delta}\Sg$ denote the group of symplectomorphisms of
$(\Sg,\om)$ equipped with the discrete topology.
Then, for any $k\geq 1$, the power
$
e_1^k\in H^{2k}(\BS^{\delta}\Sg;\bQ)
$
of the first Mumford--Morita--Miller class $e_1$ is non-trivial
for all $g\geq 3k$.
\end{theorem}
Thus, we are left with the following open problem.
\begin{problem}
Determine whether the second Mumford--Morita--Miller class
$e_2$ is non-trivial in $H^{4}(\BD_{+}^{\delta}\Sg;\bR)$
(or in $H^{4}(\BS^{\delta}\Sg;\bR)$).
\end{problem}
More generally, one can ask about polynomials in $e_{1}$ and $e_{2}$.

In the case of surfaces, a symplectic form is simply an area form and
the group of symplectomorphisms is the same as that of
area-preserving diffeomorphisms. In this paper,
we prefer to use the terminology of
symplectic topology rather than that of volume-preserving
diffeomorphisms because some of our results can be
extended to higher dimensional manifolds in the context of the
former rather than the latter.

\section{Proof of Theorem~\ref{th:main}}\label{s:proof1}

In this section we prove Theorem~\ref{th:main} by constructing
foliated surface bundles with non-zero signatures. In fact, we
prove more than was stated in Theorem~\ref{th:main}, in that we show
that {\it any} surface bundle over a surface can be made flat by
fiber summing with a trivial bundle.

First we treat the case where there is no constraint on
the total holonomy group in $\Diff_{+}\Sg$. Let
$\pi : E \ra \Sigma_h$ be any oriented $\Sg$-bundle over a closed
oriented surface of genus $h$, for example one with $\sign E\not= 0$.
Such bundles are classified by their monodromy homomorphisms
$$
\rho : \pi_1 \Sigma_h\lra \Mg \ .
$$
Choose a standard system
$\al_1,\cdots,\al_h,\beta_1,\cdots,\beta_h$ of
generators for $\pi_1\Sigma_h$ with a unique relation
$$
[\al_1,\beta_1]\cdots [\al_h,\beta_h]=1 \ ,
$$
and set
\begin{align*}
& \tilde{\al}_i = \text{any lift of}\ \rho(\al_i)\in\Mg \ \text{to}\
\Diff_{+}\Sg \\
& \tilde{\beta}_i = \text{any lift of}\ \rho(\beta_i)\in\Mg \ \text{to}\
\Diff_{+}\Sg.
\end{align*}
Then clearly we have
$$
\xi=[\tilde{\al}_1,\tilde{\beta}_1]\cdots [\tilde{\al}_h,\tilde{\beta}_h]
\in \Diff_0\Sg \ .
$$
According to a special case of a deep theorem of Thurston~\cite{Thurston},
the group $\Diff_0\Sg$ is perfect (and simple).
Hence the above element $\xi$ can be written as a product of
commutators of elements of $\Diff_0\Sg$:
$$
\xi=[\varph_1,\psi_1]\cdots [\varph_{h'},\psi_{h'}]\quad (\varph_i, \psi_i\in
\Diff_0\Sg) \ .
$$
By considering the surface $\Sigma_{h+h'}$ of genus $h+h'$ as the connected 
sum $\Sigma_h \ \sharp\ \Sigma_{h'}$, we can define a homomorphism
$$
\tilde{\rho}:\pi_1\Sigma_{h+h'}\lra \Diff_{+}\Sg
$$
by using $\tilde{\al}_i$ and $\tilde{\beta}_i$ on $\Sigma_h\setminus D^2$
and the elements $\varph_i, \psi_i\in \Diff_0\Sg$ on
$\Sigma_{h'}\setminus D^2$. Let
$$
\tilde{\pi}:\tilde{E}\lra \Sigma_{h+h'}
$$
be the corresponding foliated $\Sg$-bundle. Topologically,
if we ignore the horizontal foliation, this new bundle is
nothing but the fiber sum of the original bundle and the product
$\Sg$-bundle $\Sigma_{h'}\times\Sg$.
Hence, by Novikov additivity we have
$$
\sign \tilde{E}=\sign E\not= 0 \ .
$$
This proves Theorem \ref{th:main} in the case when we do not impose
any constraint on the holonomy of the horizontal foliation.

Next we prove that, in the above construction, we can
replace the group $\Diff_{+}\Sg$ by the subgroup $\Symp\Sg$ with
respect to a symplectic or area form $\omega$ on $\Sg$.
It is elementary to see that any Dehn twist on $\Sg$ can be
represented by an area-preserving diffeomorphism.
Since the mapping class group is generated by Dehn twists,
it follows that the natural map $\Symp\Sg\ra\Mg$ is surjective.
In fact, Moser's celebrated result~\cite{Moser} on isotopy of
volume-preserving diffeomorphisms implies
the stronger assertion that the inclusion
$$
\Symp\Sg\subset \Diff_{+}\Sg
$$
is a weak homotopy equivalence.
It follows that $\Symp\Sg\cap\Diff_{+}\Sg=\Symp_0\Sg$, and
we have an extension
$$
1\lra\Symp_0\Sg\lra\Symp\Sg\lra\Mg\lra 1 \ .
$$

\begin{remark}
Earle and Eells~\cite{EE} proved that $\Diff_{0}\Sg$ is contractible
for any $g\geq 2$. Hence $\Symp_0\Sg$ is also contractible by Moser's
result mentioned above, and we have isomorphisms
$$
H^*(\BS\Sg)\cong H^*(\BD_{+}\Sg)\cong H^*(\Mg) \ .
$$
Thus there is no difference between the characteristic classes of
smooth surface bundles and those of symplectic surface bundles,
and they are all detected by the cohomology of the mapping class group.
However, if we endow the groups
$\Diff_{+}\Sg$ and $\Symp\Sg$ with the {\it discrete} topology, then
the situation is completely different. This is the main concern of the
present paper.
\end{remark}

Now going back to the construction above,
we replace $\Diff_{+}\Sg$ by $\Symp\Sg$ and set
\begin{align*}
& \tilde{\al}_i = \text{any lift of}\ \rho(\al_i)\in\Mg \ \text{to}\
\Symp\Sg \\
& \tilde{\beta}_i = \text{any lift of}\ \rho(\beta_i)\in\Mg \ \text{to}\
\Symp\Sg \ .
\end{align*}
Then the element
$$
\xi=[\tilde{\al}_1,\tilde{\beta}_1]\cdots [\tilde{\al}_h,\tilde{\beta}_h]
$$
belongs to $\Symp_0\Sg$, not just to $\Diff_0\Sg$. But now, the group
$\Symp_0\Sg$ is not perfect. In fact, it is known that there is a
surjective homomorphism
$$
\Flux:\Symp_0\Sg\lra H^1(\Sg;\bR) \ ,
$$
called the flux homomorphism, whose kernel is the subgroup $\Ham\Sg$
consisting of Hamiltonian symplectomorphisms of $\Sg$. Fortunately
$\Ham\Sg$ is known to be perfect by a general result of 
Thurston~\cite{Thurston0}
on the group of volume-preserving diffeomorphisms of closed manifolds
(which was generalized to the case of
closed symplectic manifolds by a theorem of Banyaga~\cite{Banyaga}.
See also the books~\cite{B,MS} for these well-known results.)
Thus, we have an extension
$$
1\lra\Ham\Sg\lra\Symp_0\Sg\ \overset{\Flux}{\lra}\
H^1(\Sg;\bR)\lra 1 \ .
$$

In our situation, if $\Flux(\xi)=0$, then
$\xi$ belongs to the perfect group $\Ham\Sg$ and we are done.
In general, we cannot expect this and we have to kill
$\Flux(\xi)\in H^1(\Sg;\bR)$ in some way.
Since $\Ham\Sg$ is perfect, it is easy to see that the flux
homomorphism gives an isomorphism
$$
H_1(\Symp_0^{\delta}\Sg;\bZ)\cong H^1(\Sg;\bR) \ .
$$
The natural action of $\Symp\Sg$ on $H_1(\Symp_0^{\delta}\Sg;\bZ)$
by outer conjugation factors through that of the mapping class
group $\Mg$ because any inner automorphism of a group acts trivially
on its integral first homology, i.~e.~its abelianization.
\begin{lemma}
The flux homomorphism $\Flux:\Symp_0\Sg\ra H^1(\Sg;\bR)$ is
equivariant with respect to the natural actions of $\Mg$.
In other words, for any two elements
$\varph\in\Symp\Sg$ and $\psi\in\Symp_0\Sg$, we have the identity
$$
\Flux(\varph\psi\varph\inv)=\bar{\varph}(\Flux(\psi))
$$
where $\bar{\varph}\in\Mg$ denotes the mapping class of $\varph$ and
$\Mg$ acts on $H^1(\Sg;\bR)$ from the left by the rule
$\bar{\varph}(w)=(\bar{\varph}^{-1})^*(w)\ (w\in H^1(\Sg;\bR))$.
\label{lemma:flux}
\end{lemma}
\begin{proof}
Recall that the flux homomorphism (for the case
$\Sg\ (g\geq 2)$) can be defined as follows.
For any element $\psi\in\Symp_0\Sg$, choose an isotopy
$\psi_t\in\Symp_0\Sg$ such that $\psi_0=\id$ and
$\psi_1=\varph$. Then
$$
\Flux (\psi)=\int_{0}^{1} i_{\dot{\psi}_t}\om \ dt\in H^1(\Sg;\bR) \ .
$$
The assertion follows easily from this.
\end{proof}

As usual, let $H^1(\Sg;\bR)_{\Mg}$ denote the group of co-invariants
of $H^1(\Sg;\bR)$ with respect to the action of $\Mg$. This is
the quotient of $H^1(\Sg;\bR)$ by the subgroup generated by the elements
of the form $\varph(w)- w\ (\varph\in\Mg, w\in H^1(\Sg;\bR))$.
Notice that we have to consider
$H^1(\Sg;\bR)$ as an abelian group rather than a vector space
so that the action of $\Mg$ on it is far from being irreducible.
However, we have the following simple fact.
\begin{lemma}\label{lemma:coinv}
For any $g\geq 1$, we have $H^1(\Sg;\bR)_{\Mg}=0$.
\end{lemma}
\begin{proof}
Let $u\in H_1(\Sg;\bZ)$ be the homology class represented by
any oriented non-separating simple closed curve on $\Sg$.
Then it is easy to see that there exist elements $\varph\in\Mg$ and
$v\in H_1(\Sg;\bZ)$ such that $u=\varph(v) - v$ (consider the Dehn
twist along a non-separating simple closed curve
which intersects $u$ transversely and at only one point).
The assertion follows easily from this fact.
Moreover, it can be shown that any element in $H^1(\Sg;\bR)$
can be represented as the sum of at most $2g$ elements of the
form $\varph(w)-w$.
\end{proof}

With the above preparation, we can now finish the proof of
Theorem~\ref{th:main}.
By Lemma~\ref{lemma:coinv}, there exist elements
$\varph_i\in\Mg, w_i\in H^1(\Sg;\bR)\ (1\leq i\leq 2g)$ such that
\begin{equation}\label{eqn:flux}
\Flux(\xi)=\sum_{i=1}^{2g} (\varph_i(w_i)-w_i) \ .
\end{equation}
On the other hand, since the flux homomorphism is surjective,
for any $i$ there exists an element $\psi_i\in\Symp_0\Sg$ such that
$\Flux(\psi_i)=w_i$. By Lemma~\ref{lemma:flux}
$$
\Flux(\tilde{\varph}_i\psi_i\tilde{\varph}_i\inv)=\varph_i(\Flux(\psi_i))=
\varph_i(w_i) \ ,
$$
where $\tilde{\varph}_i\in\Symp\Sg$ is any lift of $\varph_i$.
Since $\Flux$ is a homomorphism, we can conclude
\begin{equation}\label{eqn:flux2}
\Flux([\tilde{\varph}_i,\psi_i])=\Flux(\tilde{\varph}_i
\psi_i\tilde{\varph}_i\inv)
+\Flux(\psi_i\inv)=\varph_i(w_i)-w_i \ .
\end{equation}
Now consider the element
\begin{equation}\label{eqn:eta}
\eta=[\tilde{\varph}_1,\psi_1]\cdots
[\tilde{\varph}_{2g},\psi_{2g}]\in\Symp_0\Sg \ .
\end{equation}
It follows from the equalities~\eqref{eqn:flux} and~\eqref{eqn:flux2} that
$$
\Flux(\eta)=\Flux(\xi) \ .
$$
Hence $\Flux(\xi\eta\inv)=0$ so that we have $\xi\eta\inv\in\Ham\Sg$.
Since $\Ham\Sg$ is perfect, $\xi\eta\inv$ can be represented as a
product of commutators of elements of $\Ham\Sg$. Let $h'$ be the 
number of commutators needed for this.
Then a similar argument as before yields a homomorphism
$$
\pi_1\Sigma_{h+2g+h'}\lra \Symp\Sg
$$
such that the corresponding foliated $\Sg$-bundle
$\tilde{\pi}:\tilde{E}\ra \Sigma_{h+2g+h'}$ has
total holonomy group in $\Symp\Sg$.
Now the part of $\tilde E$ over $\Sigma_h$ is the same as $E$ and
the part of $\tilde E$ over $\Sigma_{h'}$ is topologically
trivial. The remaining part of $\tilde E$ over
$\Sigma_{2g}$ may be non-trivial topologically.
However, its monodromy homomorphism to the mapping class group
factors through a free group because the mapping class of
$\psi_i$ is trivial for any $i$, and so
its signature vanishes.
Hence Novikov additivity implies that
$\sign \tilde{E}=\sign E\not= 0$.
This completes the proof of Theorem~\ref{th:main}.

\subsection{Interpretation of Theorem~\ref{th:main} in terms
of group homology}\label{ss:grouphom}
Theorem~\ref{th:main} can be translated into algebraic terms in the
context of group homology.
The extension
$$
1\lra \Diff_0\Sg\lra\Diff_{+}\Sg\lra\Mg\lra 1
$$
gives rise to the $5$-term exact sequence
\begin{align*}
H_2(\Diff_{+}^{\delta}\Sg) \lra H_2(\Mg)\lra
& H_1(\Diff_0^{\delta}\Sg)_{\Mg}\lra \\
& H_1(\Diff_{+}^{\delta}\Sg)\lra H_1(\Mg)\lra 0
\end{align*}
of integral homology groups of discrete groups.
 From Thurston's theorem~\cite{Thurston} that $\Diff_0 M$ is perfect
for any closed manifold $M$, we see that $H_1(\Diff_0^{\delta}\Sg)=0$.
Therefore, the exact sequence implies two things. Firstly, for all
$g\geq 3$, the group $\Diff_{+}\Sg$ is perfect. Secondly, the map
$H_2(\Diff_{+}^{\delta}\Sg) \lra H_2(\Mg)$ is surjective. We know from 
work of Harer that
$H_2(\Mg)\cong \bZ$ for any $g\geq 4$ and that the generator is detected
by the first Mumford--Morita--Miller class $e_1\in H^2(\Mg;\bZ)$,
see~\cite{Meyer,Harer1,Harer1a,KS}. This also holds for $g=3$, except that
$H_2(\M_{3})$ may have an additional torsion summand.
Hence we conclude that the homomorphism
$$
H_2(\Diff_{+}^{\delta}\Sg)\lra \bZ
$$
given by the cap product with $e_1$ is non-trivial for
any $g\geq 3$. This is equivalent to the existence of foliated
$\Sg$-bundles with non-zero signatures.

Next consider the extension
$$
1\lra \Symp_0\Sg\lra\Symp\Sg\lra\Mg\lra 1
$$
and the associated $5$-term exact sequence
\begin{align*}
H_2(\Symp^{\delta}\Sg) \lra H_2(\Mg)\lra
& H_1(\Symp_0^{\delta}\Sg)_{\Mg}\lra \\
& H_1(\Symp^{\delta}\Sg)\lra H_1(\Mg)\lra 0 \ .
\end{align*}
As mentioned above, the flux homomorphism yields an isomorphism
$H_1(\Symp_0^{\delta}\Sg)\cong H^1(\Sg;\bR)$.
Hence Lemma~\ref{lemma:coinv} implies that
$H_1(\Symp_0^{\delta}\Sg)_{\Mg}$ vanishes.
We can now conclude that the homomorphism
$$H_2(\Symp^{\delta}\Sg) \lra H_2(\Mg)$$
is surjective and that there is an isomorphism
$H_1(\Symp^{\delta}\Sg)\cong H_1(\Mg)$.
The former fact is equivalent to the existence of foliated
$\Sg$-bundles with area-preserving total holonomy
and with non-zero signatures as in Theorem~\ref{th:main}.
The latter fact implies that the natural projection
$\Symp\Sg\ra\Mg$ induces an isomorphism on the first integral
homology groups. In particular, the group $\Symp\Sg$ is
perfect for all $g\geq 3$.

\section{The transverse symplectic class and the flux
homomorphism}\label{s:proof2}

In this section we prove Theorem~\ref{th:main2}.
In particular, we show that the flux homomorphism
$$
\Flux: \Symp_0\Sg \lra H^1(\Sg;\bR)
$$
can be extended to a crossed homomorphism
$$
\eFlux: \Symp\Sg \lra H^1(\Sg;\bR)
$$
in an essentially unique way.

If we consider the flux homomorphism as an element
of $H^1(\Symp_0^\delta\Sg;H^1(\Sg;\bR))$, then Lemma~\ref{lemma:flux}
implies that it is invariant under the
canonical actions of $\Mg$. In other words, we can write
\begin{equation}
\Flux\in H^1(\Symp_0^{\delta}\Sg;H^1(\Sg;\bR))^{\Mg} \ .
\label{eqn:fluxm}
\end{equation}

Now we consider the cohomology class
$e+v\in H^2(\ES^{\delta}\Sg;\bR)$ mentioned in Section~\ref{s:intro}.
As we noted there, this class restricts to
$0$ on the fiber of the extension
$\pi_1\Sg\ra\ES^{\delta}\Sg\ra\Symp^{\delta}\Sg$.
Hence, in the spectral sequence $\{E^{p,q}_r\}$
for its real cohomology, we have the natural projection
\begin{align*}
p:\Ker&\left(H^2(\ES^{\delta}\Sg;\bR)\ra
H^2(\Sg;\bR)\right)
\ni e+v\\
&\lra p(e+v)\in E_{\infty}^{1,1}
\subset H^1(\Symp^{\delta}\Sg;H^1(\Sg;\bR)) \ .
\end{align*}
To prove Theorem~\ref{th:main2}, we first show the
following: if we pull back $p(e+v)$
to
$$
H^1(\Symp^{\delta}_0\Sg;H^1(\Sg;\bR))\cong
\Hom(\Symp_0\Sg,H^1(\Sg;\bR)) \ ,
$$
then we have the equality
\begin{equation}
p(e+v)=\Flux:\Symp_0\Sg\lra H^1(\Sg;\bR) \ .
\label{eqn:df}
\end{equation}
To see this, it suffices to prove the following:
\begin{lemma}\label{lemma:v}
Let $I=[0,1]$. For any $\varph\in\Symp_0\Sg$ let
$\pi: M_{\varph}\ra S^1$ be the foliated $\Sg$-bundle
over $S^1$ with monodromy $\varph$. It is the quotient
space of $\Sg\times I$ by the equivalence relation
$(p,0)\sim (\varph(p),1)$. By assumption,
there is an isotopy $\varph_{t}\in \Symp_0\Sg$
such that $\varph_{0}=\id$ and $\varph_{1}=\varph$.
Let $f : M_{\varph}\rightarrow \Sg\times S^1$ be the induced
diffeomorphism given by the correspondence
$$
M_{\varph}\ni (p,t) \longmapsto (\varph_{t}^{-1}(p),t)\in \Sg\times
S^1 \ .
$$

Then the transverse symplectic class
$v\in H^2(M_{\varph};\bR)$ is equal to
\begin{align*}
(2g-2)\mu + \Flux(\varph)\otimes \nu
\in &H^2(\Sg\times S^1;\bR)\\
\cong &H^2(\Sg;\bR)\oplus
(H^1(\Sg;\bR)\otimes H^1(S^1;\bR))
\end{align*}
under the above isomorphism, where $\mu\in H^2(\Sg;\bR)$
and $\nu\in H^1(S^1;\bR)$
denote the fundamental cohomology classes of
$\Sg$ and $S^1$ respectively.
\label{lem:fluxe}
\end{lemma}
\begin{proof}
The foliation on $M_{\varph}$ is induced from the trivial
foliation $\{\Sg\times \{t\}\}$ on $\Sg\times I$.
Hence the transverse symplectic class $v$ is represented by
the form $p^* \om$ on $\Sg\times I$,
where $p:\Sg\times I\ra \Sg$ denotes the projection
to the first factor. It is clear that the $H^2(\Sg;\bR)$-component of $v$
is equal to $(2g-2)\mu$ so that we only need to prove that
for any closed oriented curve $\ga\subset \Sg$, the value of
$v$ on the cycle $f\inv(\ga\times S^1) \subset M_{\varph}$ is equal to
$\Flux(\varph)([\ga])$ where $[\ga]\in H_1(\Sg;\bZ)$ denotes
the homology class of $\ga$. Now on $\Sg\times I$, the above
cycle is expressed as the image of the map
$$
\ga\times I\ni (q,t)\longmapsto (\varph_t(q),t)\in \Sg\times I
$$
because $f\inv(q,t)=(\varph_t(q),t)\ ((q,t)\in\Sg\times S^1)$.
Hence the required value is equal to the symplectic area
of the image of the mapping
$$
\ga\times I\ni (q,t)\longmapsto \varph_t(q)\in\Sg \ .
$$
But this is exactly equal to the value of $\Flux(\varph)$ on the
homology class represented by the cycle $\ga\subset\Sg$.
This completes the proof.
\end{proof}

Now we can finish the proof of Theorem~\ref{th:main2} as follows.
The extension
$$
1\lra\Symp_0\Sg\lra\Symp\Sg\lra\Mg\lra 1
$$
gives rise to the exact sequence
\begin{equation}
\begin{split}
0 &\lra H^1(\Mg;H^1(\Sg;\bR))  \lra H^1(\Symp^{\delta}\Sg;H^1(\Sg;\bR))\\
& \lra H^1(\Symp_0^{\delta}\Sg;H^1(\Sg;\bR))^{\Mg}
\lra H^2(\Mg;H^1(\Sg;\bR))\lra
\end{split}
\label{eqn:hh}
\end{equation}
Equations~\eqref{eqn:fluxm} and~\eqref{eqn:df}
show that the element $p(e+v)\in H^1(\Symp^{\delta}\Sg;H^1(\Sg;\bR))$
is mapped to $\Flux\in H^1(\Symp_0^{\delta}\Sg;H^1(\Sg;\bR))^{\Mg}$
in the above sequence \eqref{eqn:hh}.
In other words, the flux homomorphism can be
lifted to a crossed homomorphism
$$
\eFlux: \Symp\Sg \lra H^1(\Sg;\bR) \ .
$$
On the other hand, it was proved in~\cite{Morita3}
that $H^1(\Mg;H^1(\Sg;\bZ))=0$ for any $g\geq 1$.
It follows that $H^1(\Mg;H^1(\Sg;\bR))=0$ for any $g\geq 1$.
The exact sequence~\eqref{eqn:hh} now shows that the above
lift is essentially unique.
This completes the proof of Theorem~\ref{th:main2}.

The cohomology group $H^1(\Symp^{\delta}\Sg;H^1(\Sg;\bR))$ can be
completely determined as follows.

\begin{proposition}\label{prop:h1}
For any $g\geq 2$, there exists an isomorphism
$$
H^1(\Symp^{\delta}\Sg;H^1(\Sg;\bR))\cong \Hom_{\bQ}(\bR,\bR)
$$
where the right-hand side denotes the $\bQ$-vector space
consisting of all $\bQ$-linear mappings $\bR\ra \bR$.
Under this isomorphism, the element
$$
[\widetilde{\Flux}]\in H^1(\Symp^{\delta}\Sg;H^1(\Sg;\bR))
$$
corresponds to $\id\in \Hom_{\bQ}(\bR,\bR)$.
\end{proposition}
\begin{proof}
Consider the exact sequence~\eqref{eqn:hh}.
As mentioned above, we know that
$$
H^1(\Mg;H^1(\Sg;\bR))=0\quad (g\geq 1).
$$
On the other hand, we have the vanishing result
$$
H^2(\Mg;H^1(\Sg;\bQ))=0\quad
(g\geq 1, g\not=4, 5) \ .
$$
This is a special case of a general result of
Looijenga~\cite{Looijenga2} (for a stable range $g\geq 6$),
while the case $g\geq 9$ was already mentioned in
\cite{Morita3b}.
(See also Proposition~21 of \cite{HR}. The proof there should be
modified to use Harer's result~\cite{Harer2} on the {\it third} homology
group of the moduli spaces as well as results of
Igusa~\cite{Igusa} and Looijenga~\cite{Looijenga1} for
low genera $g=2, 3$, instead of Harer's earlier result~\cite{Harer1}
on the second homology. This correction forces us to exclude $g=4$ or 
$5$ for now.)

Thus we have an isomorphism
\begin{equation}
H^1(\Symp^{\delta}\Sg;H^1(\Sg;\bR))\cong
H^1(\Symp_0^{\delta}\Sg;H^1(\Sg;\bR))^{\Mg} \ ,
\label{eqn:hh1}
\end{equation}
except possibly for $g=4, 5$ for the moment.

Now the flux homomorphism gives rise to an isomorphism
$$
\Flux: H_1(\Symp_0^{\delta}\Sg;\bZ)\cong H^1(\Sg;\bR) \ .
$$
Hence we can write
$$
H^1(\Symp_0^{\delta}\Sg;H^1(\Sg;\bR))\cong
\Hom(H^1(\Sg;\bR),H^1(\Sg;\bR))
$$
and under this isomorphism, the flux homomorphism clearly
corresponds to the identity.
On the other hand, an analysis of the action of $\Mg$ on the
right-hand side yields an isomorphism
\begin{equation}
\Hom(H^1(\Sg;\bR),H^1(\Sg;\bR))^{\Mg}\cong \Hom_{\bQ}(\bR,\bR) \ .
\label{eqn:Hamel}
\end{equation}
More precisely, if we choose a Hamel basis
$\{a_{\la}\}_{\la}$ for $\bR$ considered as a vector space
over $\bQ$, then we have an isomorphism
$$
H^1(\Sg;\bR)\cong \oplus_{\la} H^1(\Sg;a_{\la}\bQ) \ .
$$
It is easy to see that any endomorphism
$f:H^1(\Sg;\bR)\ra H^1(\Sg;\bR)$, which is equivariant with
respect to the natural action of $\Mg$ must send any
summand $H^1(\Sg;a_{\la}\bQ)$ to a direct sum of finitely
many such summands by some scaler multiplication in each factor.
The isomorphism~\eqref{eqn:Hamel} follows from this.

We can eliminate the possible exceptions for~\eqref{eqn:hh1} by a 
stabilization argument using the simple behavior of the flux 
homomorphisms under the inclusions
$$
\Symp_0^c \Sigma_g^0\subset \Symp_0\Sg \ ,
\quad
\Symp_0^c \Sigma_g^0\subset\Symp_0^c\Sigma_{g+1}^0 \ .
$$
Here $\Sigma_g^0=\Sg\setminus D^2$ and
$\Symp_0^c \Sigma_g^0$ denotes the group of symplectomorphisms
of $\Sigma_g^0$ with compact supports (see the next section
for the flux homomorphism for the group $\Symp_0^c \Sigma_g^0$).

Thus the isomorphism~\eqref{eqn:hh1}
holds for any $g\geq 2$
and the required result follows.
\end{proof}

\begin{remark}
Kawazumi kindly pointed out the following simple argument which avoids
the use of the vanishing result for $H^2(\Mg;H^1(\Sg;\bQ))$.
Any element in
$$
H^1(\Symp_0^{\delta}\Sg;H^1(\Sg;\bR))^{\Mg}\cong
\Hom(H^1(\Sg;\bR),H^1(\Sg;\bR))^{\Mg}
$$
is obtained from the flux homomorphism
(= the identity) by applying some
endomorphism to the coefficients $H^1(\Sg;\bR)$ which is
equivariant under the action of $\Mg$.
Since we already proved that the identity can be lifted,
any other element can also be lifted
simply by applying an $\Mg$-equivariant change of coefficients.
\end{remark}

\section{The symplectomorphism groups of open surfaces}\label{s:open}

In this section, we prepare a few facts concerning the
symplectomorphism groups of {\it open} surfaces.
These will be used in the proof of Theorem~\ref{th:main3}
given in the next section.

Let $D^2\subset\Sg$ be a closed embedded disk and $\Sigma_g^{0}$
the open surface $\Sg\setminus D^2$.
We consider the group $\Symp^{c}\Sigma_g^{0}$ of
symplectomorphisms of $\Sigma_g^{0}$ with {\it compact supports}.
Let $\Symp^{c}_{0}\Sigma_g^{0}$ denote the identity component of
$\Symp^{c}\Sigma_g^{0}$. In this context,
we again have a flux homomorphism
$$
\Flux: \Symp^{c}_{0}\Sigma_g^{0}\lra H^1_{c}(\Sigma_g^{0};\bR) \ ,
$$
where $H^1_{c}(\Sigma_g^{0};\bR)$ denotes the first real cohomology
group of $\Sigma_g^{0}$ with compact supports.
It is easy to see that the inclusion $\Sigma_g^{0}\subset \Sg$
induces an isomorphism $H^1_{c}(\Sigma_g^{0};\bR)\cong H^1(\Sg;\bR)$,
and that the following diagram is commutative:
\begin{equation}
\begin{CD}
\Symp^{c}_{0}\Sigma_g^{0} @>{\Flux}>> H^1_{c}(\Sigma_g^{0};\bR)\\
@VVV @VV{\cong}V\\
\Symp_{0}\Sigma_g @>>{\Flux}> H^1(\Sg;\bR).
\end{CD}
\label{eqn:FFlux}
\end{equation}
It is known that the kernel of the flux homomorphism is
equal to the subgroup $\Ham^{c}\Sigma_g^{0}$ consisting of
Hamiltonian symplectomorphisms with compact supports.
Thus we have an extension
$$
1\lra\Ham^{c}\Sigma_g^{0}\lra\Symp^{c}_{0}\Sigma_g^{0}\
\overset{\Flux}{\lra}\ H^1_{c}(\Sigma_g^{0};\bR)\lra 1 \ .
$$
In contrast to the case of closed surfaces, the group
$\Ham^{c}\Sigma_g^{0}$ is not perfect. In fact,
there is a surjective homomorphism
$$
\CAL:\Ham^{c}\Sigma_g^{0}\lra \bR \ ,
$$
called the (second) Calabi homomorphism, see~\cite{Calabi}.
The kernel of this homomorphism is known to be
simple, and hence perfect, by a result of Banyaga~\cite{Banyaga,B}.

The flux and the Calabi homomorphisms can be defined for any
non-compact symplectic manifold. Here we only consider the case
of {\it exact} symplectic manifolds, assuming the
existence of a $1$-form $\la$ such that $\om=-d\la$.
The open surface $\Sigma_g^{0}$ which we are concerned with
is an exact symplectic manifold.

For any exact symplectic manifold $(M,\om)$ of dimension
$2n$, the flux and the Calabi homomorphisms can be expressed as
$$
\Flux(\varph)=[\la-\varph^*\la]\in H^1_{c}(M;\bR)\quad
(\varph\in\Symp_0^{c} M)
$$
and
\begin{equation}
\CAL(\varph)=-\frac{1}{n+1}\int_{M} \varph^*\la\land \la\land
\om^{n-1}\quad (\varph\in \Ham^{c} M)
\label{eqn:CAL}
\end{equation}
respectively (see Lemma~10.14 and Lemma~10.27 of~\cite{MS}).
The formula~\eqref{eqn:CAL} above can be used for any
$\varph\in\Symp_0^{c}M$, not necessarily in $\Ham^{c} M$. It defines a
{\it map}
\begin{equation}
\CAL:\Symp_0^{c} M\lra \bR \ ,
\label{eqn:CAL2}
\end{equation}
and a straightforward calculation shows that
\begin{equation*}
\CAL(\varph\psi)=\CAL(\varph)+\CAL(\psi)
+\frac{1}{n+1}\int_{M} \Flux(\varph)\land \Flux(\psi)\land \om^{n-1}
\end{equation*}
for any two elements $\varph,\psi\in \Symp_0^{c} M$.
Hence the map~\eqref{eqn:CAL2} above is
a homomorphism if and only if the pairing
\begin{equation}
H^1_{c}(M;\bR)\otimes H^1_{c}(M;\bR)\ni
([\alpha],[\beta])\longmapsto \int_{M} \alpha\land\beta\land \om^{n-1}
\label{eqn:pairing}
\end{equation}
is trivial. This is the case if $\dim M=2n\geq 4$,
because then the integrand is exact (and compactly supported).
However, in our case of an open surface $M=\Sigma_g^{0}$,
the pairing~\eqref{eqn:pairing} is non-trivial, and even non-degenerate.
Define the Heisenberg group $\mathcal{H}$ to be the central extension
$$
0\longrightarrow\bR\longrightarrow\mathcal{H}\longrightarrow
H^1_{c}(\Sigma_g^{0};\bR)\longrightarrow 1
$$
corresponding to the cup product pairing
$H^1_{c}(\Sigma_g^{0};\bR)\otimes H^1_{c}(\Sigma_g^{0};\bR)\ra\bR$.
Now we obtain the following fact:
\begin{proposition}\label{prop:os}
For any $g\geq 2$, the mapping $\CAL+\Flux$ defines a
surjective homomorphism
$$
\CAL+\Flux: \Symp_0^{c}\Sigma_g^{0}\ \lra\ \mathcal{H} \ .
$$
\end{proposition}

\begin{corollary}\label{cor:h1}
The flux homomorphism induces an isomorphism
$$
\Flux: H_1(\Symp^{c}_{0}\Sigma_g^{0};\bZ) \cong
H^1_{c}(\Sigma_g^{0};\bR) \ .
$$
Furthermore for any real number $r\in \bR$, there exist
two elements $\varph, \psi\in \Symp^c_{0}\Sigma_g^0$
such that the commutator $\eta=[\varph,\psi]\in\Ham^c \Sigma_g^0$
satisfies $\CAL(\eta)=r$.
\end{corollary}
\begin{proof}
The first statement follows from Proposition~\ref{prop:os} together with
the fact that $\Ker\CAL\subset Ham^{c}\Sigma_g^{0}$
is perfect. The second statement follows easily from the above argument.
\end{proof}

\section{Proof of Theorem~\ref{th:main3}}\label{s:proof3}

In this section we prove Theorem~\ref{th:main3}, which
shows the non-triviality of any power
$e_1^k \in H^{2k}(\BS^{\delta}\Sg;\bR)$ of the first
Mumford--Morita--Miller class whenever $g\geq 3k$.

We first treat the case where the total holonomy group is in 
$\Diff_{+}\Sg$, rather than in $\Symp\Sg$. As in the previous section,
fix a closed embedded disk $D^2\subset\Sg$.
We denote by $\Diff (\Sg,D^2)$ the group of diffeomorphisms
of $\Sg$ which are the identity on some open neighborhoods
of $D^2$. This is the same as the group $\Diff^{c}\Sigma_g^0$
of diffeomorphisms with compact supports
of the open surface $\Sigma_g^0$.

Let $\pi:E\ra \Sigma_h$ be any $\Sg$-bundle over $\Sigma_h$, for
example one with $\sign E\not= 0$. Then we can apply the same
argument as in Section~\ref{s:proof1} to this bundle
replacing the group $\Diff_{+}\Sg$ by $\Diff (\Sg,D^2)$.
Fortunately Thurston's theorem (see~\cite{Thurston,B})
is also valid for this relative case, giving that the identity
component $\Diff_0 (\Sg,D^2)$ is simple and hence perfect.
Thus for some $h'$, we obtain a homomorphism
$$
\pi_1\Sigma_{h+h'}\lra \Diff (\Sg,D^2)
$$
such that the signature of the total space of the
associated foliated $\Sg$-bundle over $\Sigma_{h+h'}$
is equal to $\sign E\not= 0$.
This implies the non-triviality
$$
e_1\not= 0\in H^2(\BD^{\delta}(\Sg,D^2);\bQ) \ .
$$

To prove the non-triviality of higher powers $e_1^k$,
consider a genus $kg$ surface
$\Sigma_{kg,1}=\Sigma_{kg}\setminus \Int D^2$
with one boundary component as the boundary connected sum
$$
\Sigma_{kg,1}=\Sigma_{g,1}\ \natural \cdots\natural\ \Sigma_{g,1}
$$
of $k$ copies of $\Sigma_{g,1}=\Sg\setminus \Int D^2$.
This induces a homomorphism
\begin{equation}
f_k: \Diff (\Sg,D^2)\times\cdots\times \Diff (\Sg,D^2)
\lra \Diff (\Sigma_{kg},D^2)
\label{eqn:fk}
\end{equation}
from the direct product of $k$ copies of the group
$\Diff (\Sg,D^2)$ to $\Diff (\Sigma_{kg},D^2)$.
It can be shown that
$$
f_{k}^*(e_1)=e_1\times 1\times\cdots\times 1+\cdots
+1\times\cdots\times 1\times e_1 \ ,
$$
see~\cite{Miller,Morita2} or~\cite{Morita4}.
It follows that
$$
f_{k}^*(e_1^k)=e_1\times\cdots\times e_1+\text{other terms}
$$
where the other terms belong to various summands of
$$
H^*(\BD^{\delta}(\Sg,D^2);\bQ)\otimes\cdots\otimes 
H^*(\BD^{\delta}(\Sg,D^2);\bQ)
$$
other than
$$
H^2(\BD^{\delta}(\Sg,D^2);\bQ)\otimes\cdots\otimes H^2(\BD^{\delta}(\Sg,D^2);\bQ).$$Since $e_1\times\cdots\times e_1\not= 0$, we 
can now conclude that
$f_{k}^*(e_1^k)\not= 0$.
This proves the non-triviality
$$
e_1^k\not= 0\in H^{2k}(\BD_{+}\Sg;\bQ)\quad
\text{for any $g\geq 3k$} \ .
$$

Next we consider the case where the total holonomy is contained
in $\Symp\Sg$. We apply the same argument as in Section~\ref{s:proof1},
but replacing the group $\Symp\Sg$ by $\Symp^c\Sigma_g^0$.
At the final stage, we must use the second statement of
Corollary~\ref{cor:h1} to kill the value of the Calabi homomorphism.
To summarize, we kill the value of the flux homomorphism by adding
$2g$ commutators in $\Symp^c\Sigma_g^0$ as in Section~\ref{s:proof1} and then
kill the value of the Calabi homomorphism by adding one commutator
in $\Symp^c_{0}\Sigma_g^0$. Then we can use the perfection
of the subgroup $\Ker\CAL\subset \Symp^c_{0}\Sigma_g^0$
to show the non-triviality
$$
e_1\not= 0\in H^2(\BS^{c,\delta}\Sigma_g^0;\bQ) \ .
$$
Finally we consider the homomorphism
$$
h_k: \Symp^{c}\Sigma_g^0\times\cdots\times \Symp^{c}\Sigma_g^0
\lra \Symp^{c}\Sigma_{kg}^0
$$
which is defined similarly to the $f_k$ in~\eqref{eqn:fk} and
apply the same argument as above to show
the non-triviality of the power $e_1^k$.

This completes the proof of Theorem~\ref{th:main3}.

\section{Further results}\label{s:final}

\subsection{Perfect versus uniformly perfect groups}
Combining our discussion in~\ref{ss:grouphom} with the main result
of~\cite{EK}, we obtain the following:
\begin{corollary}
Let $G=\Diff_{+}\Sg$ or $\Symp\Sg$. For all $g\geq 3$ the group $G$
is perfect but not uniformly perfect.

Moreover, if $\varphi\in G$ represents the Dehn twist along any
homotopically non-trivial simple closed curve on $\Sg$, then the
commutator length of $\varphi^{k}$ in $G$ grows linearly with $k$, for
all $g\geq 2$.
\end{corollary}
\begin{proof}
     We saw in~\ref{ss:grouphom} that $H_{1}(G^{\delta})=H_{1}(\Mg)$ for all
     $g\geq 2$. For $g\geq 3$, the mapping class group is known to be
     perfect, see for example~\cite{Harer1a}.
 
     The projection $G\rightarrow\Mg$ is surjective, and the
     commutator length of $\varphi^{k}$ is bounded below by that of
     its image in $\Mg$. Thus the main result of~\cite{EK} gives the
     conclusion, compare also~\cite{BK}.
     \end{proof}
For a perfect group $G$ not being uniformly perfect is equivalent to the 
statement that the comparison map 
$c\colon H^{2}_{b}(G^{\delta})\rightarrow H^{2}(G^{\delta})$ 
from the second bounded cohomology to the usual cohomology with real 
coefficients is not injective. If we denote the kernel of $c$ by 
$K(G^{\delta})$, it is easy to see that $K(\Mg)$ injects into 
$K(G^{\delta})$ for 
$G=\Diff_{+}\Sg$ or $\Symp\Sg$. The result of~\cite{EK} to the effect 
that $K(\Mg)$ is non-zero has been generalized by Bestvina and 
Fujiwara~\cite{BF} to show that it is infinite-dimensional. Thus, 
$K(G^{\delta})$ is also infinite-dimensional.

Note that because the Mumford--Morita--Miller class $e_{1}\in 
H^{2}(\Mg)$ is a bounded class, i.~e.~is in the image of $c$, the 
same is true for $e_{1}\in H^{2}(G^{\delta})$ and its powers $e_{1}^{k}$. Thus 
Theorem~\ref{th:main3} shows in particular that 
the comparison map $c$ is non-trivial on $H^{2k}_{b}(G^{\delta})$ for 
$g\geq 3k\geq 3$ and $G=\Diff_{+}\Sg$ or $\Symp\Sg$.

The groups $\Diff_0\Sg$, $\Diff_0 (\Sg,D^2)$, $\Ham\Sg$ and
$\Ker\CAL\subset \Ham^c\Sigma_g^0$ are perfect by the results of
Thurston~\cite{Thurston0, Thurston} and Banyaga~\cite{Banyaga}, compare
also~\cite{B}. In parallel to our work on this paper, Gambaudo and
Ghys~\cite{GG} have proved that $\Ham\Sg$ is not uniformly perfect.
Their arguments also apply to the group $\Ker\CAL\subset
\Ham^c\Sigma_g^0$, although they do not state this in~\cite{GG}.
The result of Gambaudo--Ghys implies that in our proof of
Theorem~\ref{th:main} one cannot control the base genus of the trivial
fibration which we fiber sum to a given surface bundle to obtain a
flat bundle with total holonomy in $\Symp\Sg$.

Note that Entov and Polterovich~\cite{EP} recently proved that $\Ham M$ is 
not uniformly perfect if $M$ belongs to a certain subclass of the 
spherically monotone symplectic manifolds which includes $S^{2}$ and 
many high-dimensional manifolds, but not the surfaces of positive genus.

Whether or not $\Diff_0\Sg$ and $\Diff_0 (\Sg,D^2)$ are uniformly
perfect remains a very interesting open question.
 
\subsection{Symplectic pairs}
A {\it symplectic pair} on a smooth manifold is a pair of closed
two-forms $\omega_{1}$, $\omega_{2}$ of constant and complementary
ranks, for which $\omega_{1}$ restricts as a symplectic form to the
leaves of the kernel foliation of $\omega_{2}$, and vice versa. This
definition is analogous to that of contact pairs and of
contact-symplectic pairs discussed by Bande~\cite{Bande1,Bande2}.

Manifolds with symplectic pairs are always symplectic, but they satisfy
much stronger topological restrictions than general symplectic manifolds.
For example, a four-manifold with a symplectic pair admits symplectic
structures for both choices of orientation, because
$\omega_{1}+\omega_{2}$ and $\omega_{1}-\omega_{2}$ are symplectic
forms inducing opposite orientations.

Theorem~\ref{th:main} implies:
\begin{corollary}
     There exist smooth closed oriented four-manifolds of non-zero
     signature which admit symplectic pairs.
     \end{corollary}
The signature vanishes for all other four-manifolds which we know to
admit symplectic pairs.

\subsection{The crossed flux homomorphism in higher dimensions}
Let $(M,\om)$ be any closed symplectic manifold and
${\mathcal M}_{\om}$ its symplectic mapping class group
defined to be the quotient of $\Symp (M,\om)$ by its
identity component $\Symp_0 (M,\om)$, so that we have an
extension
\begin{equation}
1\lra\Symp_0 (M,\om)\lra\Symp (M,\om)
\lra{\mathcal M}_{\om}\lra 1 \ .
\label{eqn:mom}
\end{equation}
In view of Theorem~\ref{th:main2}, it appears to be an interesting problem
to determine whether the flux homomorphism
\begin{equation}
\Flux : \Symp_0 (M,\om)\lra H^1(M;\bR)/\vGa_{\om}
\label{eqn:FF}
\end{equation}
can be extended to a crossed homomorphism on the whole
group $\Symp (M,\om)$ or not.
Here $\vGa_{\om}$ denotes the flux subgroup corresponding to the flux
of non-trivial loops in $\Symp_{0} (M,\om)$. (This has to be divided
out to make the flux well-defined, see~\cite{MS}.)
The extension \eqref{eqn:mom} yields an exact sequence
\begin{align*}
0 & \lra H^1({\mathcal M}_{\om};H^1(M;\bR)/\vGa_{\om})  \lra
H^1(\Symp (M,\om);H^1(M;\bR)/\vGa_{\om})\\
& \lra H^1(\Symp_0 (M,\om);H^1(M;\bR)/\vGa_{\om})^{{\mathcal M}_{\om}}
\overset{\delta}{\lra}
H^2({\mathcal M}_{\om};H^1(M;\bR)/\vGa_{\om})\lra
\end{align*}
It is easy to generalize Lemma~\ref{lemma:flux}, which
treats the case of surfaces, to the case of
closed symplectic manifolds.
Hence we can write
$$
\Flux\in H^1(\Symp_0 (M,\om);H^1(M;\bR)/\vGa_{\om})^{{\mathcal M}_{\om}}
$$
and we may ask whether the element $\delta(\Flux)$ is trivial or not.
This is equivalent to asking whether the extension
\begin{equation}
1\lra H^1(M;\bR)/\vGa_{\om}\lra\Symp(M,\om)/\Ham(M,\omega)
\lra{\mathcal M}_{\om}\lra 1 \
\label{eqn:split}
\end{equation}
splits or not.
If this is the case, then the flux extends and
the group $H^1({\mathcal M}_{\om};H^1(M;\bR)/\vGa_{\om})$
measures the differences between the possible extensions.

As a partial answer to this problem, we have the following.
Assume that the cohomology class $[\om]\in H^2(M;\bR)$ is a 
multiple of the first Chern class $c_1(M)\in H^2(M;\bZ)$. 
(This is a variant of the monotonicity assumption.)
Then it was proved by McDuff~\cite{McDuff} and by Lupton--Oprea~\cite{LO}
that the flux subgroup $\vGa_\omega$ is trivial. We can extend 
this result as follows, thereby also reproving the triviality of 
$\vGa_\omega$ from our point of view.

\begin{proposition}
Let $(M,\om)$ be a closed symplectic manifold and assume that
the cohomology class $[\om]\in H^2(M;\bR)$ is a multiple of
the first Chern class $c_1(M)\in H^2(M;\bZ)$.
Then the flux subgroup $\vGa_\omega$ is trivial and the
flux homomorphism
$\Flux: \Symp_0 (M,\om)\lra H^1(M;\bR)$ can be canonically
extended to a crossed homomorphism
$$
\widetilde{\Flux} : \Symp (M,\om)\lra H^1(M;\bR).
$$
\label{prop:fluxhi}
\end{proposition}
\begin{proof}
We modify the argument in the proof of Theorem~\ref{th:main2}, given in 
Section~\ref{s:proof2}, as follows.
Observe first that the Euler class $e\in H^2(\ES^\delta\Sg;\bZ)$
considered there is nothing but the first Chern class of the tangent bundle
along the fibers of the universal surface bundle over $\BS^\delta\Sg$.
Let $\BS^\delta (M,\omega)$ be the classifying space of the discrete group
$\Symp^\delta (M,\omega)$ and let
$$
\pi:\ES^\delta (M,\omega)\lra \BS^\delta (M,\omega)
$$
be the universal foliated $M$-bundle over $\BS^\delta (M,\omega)$
with total holonomy group in $\Symp(M,\omega)$.
Then we have the first Chern class
$$
c_1(\xi)\in H^2(\ES^\delta (M,\omega);\bZ)
$$
where $\xi$ denotes the tangent bundle along the fibers of $\pi$.
By assumption, there exists a non-zero real number $r$ such
that $[\om]=r c_1(M)$. Now consider the cohomology class
$$
u=v- r c_1(\xi)\in H^2(\ES^\delta (M,\omega);\bR)
$$
where
$v$ denotes the transverse symplectic class represented by
the global $2$-form $\tilde\om$ on $\ES^\delta (M,\omega)$ which
restricts to $\om$ on each fiber.
The restriction of $u$ to the fiber vanishes so that,
in the spectral sequence $\{E^{p,q}_r\}$
for the real cohomology, we have 
$$
p(u)\in E_{\infty}^{1,1}
\subset H^1(\BS^{\delta}(M,\omega);H^1(M;\bR)) \ .
$$
Now we consider the composition of homomorphisms
\begin{align*}
  H^1(\BS^{\delta}&(M,\omega);H^1(M;\bR))
\lra  H^1(\BS_0^{\delta}(M,\omega);H^1(M;\bR))\\
\lra & H^1(\BS_0^{\delta}(M,\omega);H^1(M;\bR)/\vGa_\omega)\\
&\cong\Hom(\Symp_0(M,\omega),H^1(M;\bR)/\vGa_\omega)
\end{align*}
where the first homomorphism is induced by the
restriction to the subgroup $\Symp_0(M,\omega)\subset \Symp(M,\omega)$
while the second one is induced by the natural
projection $H^1(M;\bR)\ra H^1(M;\bR)/\vGa_\omega$.
Let
$$
\overline{p(u)}\in
\Hom(\Symp_0(M,\omega),H^1(M;\bR)/\vGa_\omega)
$$
be the image of $p(u)$ under the above composition.
Then we have the equality
\begin{equation}
\overline{p(u)}=\Flux:\Symp_0(M,\omega)\lra H^1(M;\bR)/\vGa_\omega\ .
\label{eqn:df2}
\end{equation}
This can be shown by suitably adapting Lemma~\ref{lem:fluxe} to the
case of a general closed symplectic manifold $M$ instead of
$\Sigma_g$. Thus we see that the flux homomorphism
can be extended to a homomorphism from $\Symp_0(M,\omega)$
to the whole of $H^1(M;\bR)$ (rather than its quotient by $\vGa_\omega$).
On the other hand, Banyaga's result~\cite{Banyaga} that $\Ker 
\Flux=\Ham(M,\omega)$
is perfect (and simple) implies that the abelianization of $\Symp_0(M,\omega)$
is equal to $H^1(M;\bR)/\vGa_\omega$. We can now conclude that
the flux subgroup $\vGa_\omega$ is trivial and further that the
flux homomorphism can be extended canonically to a crossed homomorphism
on the whole group $\Symp(M,\omega)$.
This completes the proof.
\end{proof}

\begin{example}
The above proof does not apply to the torus $T^2$ with the
standard symplectic form $\omega_0$ because the first Chern class is trivial
in this case. In fact, the flux subgroup is isomorphic to
$H^1(T^2;\bZ)$ which is non-trivial.
However the flux homomorphism does extend canonically to a
crossed homomorphism
$$
\eFlux:\Symp (T^2,\omega_0)\lra H^1(T^2;\bR)/H^1(T^2;\bZ).
$$
This is because the mapping class group ${\mathcal M}_1\cong {\rm SL}(2,\bZ)$
acts on $T^2$ linearly by symplectomorphisms and hence
the extension \eqref{eqn:split} splits canonically.

\end{example}

\medskip
\noindent
\textsc{Acknowledgements:}
The second named author would like to thank the
Mathematisches Institut der Universit\"at M\"unchen,
where the present work was done,
for its hospitality. Thanks are also due to
N.~Ka\-wa\-zu\-mi, S.~Matsumoto,
K.~Ono and T.~Tsuboi for enlightening discussions and helpful information.

\bibliographystyle{amsplain}

\end{document}